\documentclass[a4paper,11.5pt,reqno]{amsart}

\usepackage{graphicx}
\usepackage[utf8]{inputenc}
\usepackage{amsmath,amssymb,amsthm}
\usepackage{enumitem}
\usepackage{xcolor}
\usepackage{hyperref}

\usepackage{lineno}
\modulolinenumbers[3]

\theoremstyle{plain}

\theoremstyle{definition}

\theoremstyle{remark}

\newcommand{\cc}      { \text{c.c.} }                     
\newcommand{\ash}     { {\color{red} \; \mathbf{+1}} }    
\newcommand{\qsh[1]}  { {\color{blue} \mathbf{#1}} }      
\newcommand{\upth}      { ^\text{th} }                    

\numberwithin{equation}{section}

\hypersetup{colorlinks = true, linkcolor = blue, urlcolor = blue, citecolor = blue}

\begin{document}

\title[Constructing the $2^{2^5} + 1 = F_5$-gon, the first step] 
      {Constructing the $\mathbf{2^{2^5} + 1 = F_5}$-gon, the first step: \\ the $641$-gon}

\author[H. Ruhland]{Helmut Ruhland}
\address{Santa F\'{e}, La Habana, Cuba}
\email{helmut.ruhland50@web.de}



\begin{abstract}
In an article, S. Adlaj showed the construction of the $F_4 = 65\,537$-gon ($F_4$ a Fermat prime). By extending the author's method to cases where the multiplicative group is not of order $2^n$ but has an additional prime factor, I here construct the $641$-gon. $641$ is the smaller prime factor of the composite Fermat number $F_5 = 2^{2^5} + 1$. 
\end{abstract}

\date{\today}

\maketitle

\section{Introduction}

The $n\upth$ Fermat number is defined as $F_n = 2^{2^n} + 1$. For $n = 1, 2, 3, 4,$ these numbers are prime. $F_5 = 4\,294\,967\,297$ is composite and equals $641 \cdot 6\,700\,417 = (2^7 \cdot 5 + 1) \, (2^7 \cdot 52\,347 + 1)$. \\
When $F_n$ is a prime, the regular $n$-gon can be constructed using only square roots. For the regular $F_4$-gon, this was shown by Hermes \cite{Hermes}, who carried out the calculation by hand. In \cite{Adlaj} the construction was shown using a CAS. The amount of time Hermes spent on the calculation by hand differs significantly from the time writing the article \cite{Adlaj} by its author. \\

The order of the multiplicative group of the primes $P$ treated in this article has the form $2^N \cdot Q, Q$ a prime $\ge 3$. So we have $N$ positions in the sequence of field extensions where the $Q$-section can be performed. Here we consider the cases where the $Q$-section is performed first, followed by the $N -1$ quadratic extensions. \\

The subsequent quadratic field extensions are done with the scheme (slightly extended) introduced in \cite{Adlaj} for the $F_4 = 65\,537$-gon. \\

In the appendices \ref{examples_begin} to \ref{examples_end}, some examples with small primes $P$ are given. The special cases $Q = 2^n$ allow the construction of the Fermat prime cases; see appendix \ref{special_fermat} for $F_2, F_3, F_4$. \\

It is highly propable that no article constructing the $6\,700\,417$-gon (the bigger prime factor of $F_5$)
will be written, for obvious reasons. Thus, the $F_5$-gon will remain unconstructed.

\section{The smaller prime factor $P = 641$ as divisor of $F_5$, \\
			the two quintisections first}

The order of the multiplicative group modulo $P$: $\phi (P) = 2^7 \cdot Q$, where $Q$ is a prime (here $Q = 5$). First, the quintisection is done. Because of the factor $2^7$, the last 6 quadratic extensions of the field $K$ with $[K : \mathbf{Q}] = 5$, can be done in almost the same manner as in \cite{Adlaj}. 'Almost' means that a second lower subscript in the following sums has to be introduced. \\

Let $s_{q,a}^b$ denote the following sum, call the upper index $b$ the level, the additional index $q$ represents the $Q$ roots gotten by the quintisection:
\begin{align}
   s_{q,a}^b := & \sum_{k = 0}^{k = (P - 1) / (Q \, 2^{b+1}) - 1} c(q + (a + k \, 2^b) \, Q), \quad c(l) := 2 \cos \left( \frac {2 \pi \, 3^l}{P} \right) \label{sqab_definition} \\
	& 0 \le q \le Q - 1, \quad  0 \le a \le 2^b - 1  \nonumber
\end{align}

where the number 3, appearing in the definition of the algebraic integer $c(\cdot)$, might be replaced with any other generator of the multiplicative group of integers modulo the prime $P$ without altering the sum, merely reordering its summands. \\

The $Q$-section yields the initial (for the following quadratic extensions) sums $s_{0,0}^0,\dots, s_{4,0}^0$. They are defined as roots of the quintic $s^5 + s^4 - 258 s^3 - 564 s^2 + 5328 s - 5120$. The roots can be calculated by $Q$-section because the cyclic group of order $5$ is the Galois group of this quintic. The solutions can be expressed using roots of unity derived from the Lagrange resolvents (see \cite{Adlaj_RootsUnity}). Here, and in the appendices the principal values of the multivalued roots have to be taken. Here, "c.c." in \ref{sq0} means the two complex conjugate terms. With quintisection, we get for the sums:
\begin{align}
	& C_{\pm} = \sqrt[5]{\frac{15\,128 \mp 7\,750 \, \sqrt{5} - (4\,275 \pm 1\,895 \, \sqrt{5}) \, \sqrt{10 \pm 2 \, \sqrt{5}} \, i} {2 \cdot 641 ^ {3/2}}}  \quad \vert C_{\pm} \vert = 1 \nonumber \\
	& s_{q,0}^0 = (- 1 + \sqrt{641} \, ( \, \epsilon^{1-q} \, C_+ 
	                                    +   \epsilon^{-2q} \, C_- + \cc \, )) \, / \, 5 \quad \epsilon = e^{2 \pi i / 5}, \quad q = 0, \dots,  4 \label{sq0}
\end{align}

Similar to \cite{Adlaj}, lower leading indices 0 can be omitted for the sake of brevity, i.e., $s_{0,a}^b$ can be written as $s_{a}^b$, and $s_{0,0}^b$ can be written as $s^b$. 
Note, in particulaupr, that the sum $s_{0,0}^b$, which we denote with $s^b$, is the sum of all $5\upth$ powers of a primitive $P\upth$ root of unity. So we have for each integer $b, 0 \le b \le 5$, the two sums $s_{q,0}^{b+1}$ and $s_{q,2^b}^{b+1}$ are roots of the quadratic polynomial
\begin{equation}
   p_{q,0}^b (x) := x ^ 2 - s_{q,0}^b \, x +  s_{q,2^b}^{b+1} \, s_{q,0}^{b+1}   
\end{equation}

Our task is accomplished once we express the algebraic integer $s^6 = 2 \cos(2 \pi / P)$, which algebraic degree over $\mathbf{Q}$ is $2^6 = 64 \, Q$, via nested quadratic radicals, which we calculate via consecutively solving $6$ quadratic equations. \\
 
We begin by observing that each sum and product of the two $s_{q,2^b}^{b+1}$ and $s_{q,0}^{b+1}$ reside in an extension field of $\mathbf{Q}$ whose algebraic degree is $2^b\, Q$. In fact, it might be expressed as an integer linear combination of algebraic integers whose algebraic degree, over $\mathbf{Q}$ does not exceed $2^b\, Q$. \\

{\large \bf Expressing $\mathbf{N - 1 = 6}$ products as integer linear combinations of algebraic integers} \\

We successively calculate:

\begin{align}
   & s_1^1 s^1 &    = \enspace  & 6 s^0 + 6 s_{1,0}^0 + 7 s_{2,0}^0 + 7 s_{3,0}^0 + 6 s_{4,0}^0 \nonumber \\
   & s_2^2 s^2 &    = \enspace  & (3 s^1 + 2 s_1^1)  +  (2 s_{1,0}^1 + s_{1,1}^1) + (2 s_{2,0}^1 + 2 s_{2,1}^1) + (s_{3,0}^1 + s_{3,1}^1) + (2 s_{4,0}^1) \nonumber \\
   & s_4^3 s^3 &    = \enspace  & (s_{2}^2) + (s_{1,0}^2 + 2 s_{1,2}^2 + s_{1,3}^3) + (3 s_{2,0}^3) \nonumber \\
   & s_8^4 s^4 &    = \enspace  & (s_{1,6}^3) + (s_{3,2}^3 + s_{3,6}^3 + s_{3,7}^3) \\
   & s_{16}^5 s^5 & = \enspace  & (s_{1,2}^4) + (s_{4,12}^4) \nonumber \\
   & s_{32}^6 s^6 & = \enspace  & s_{31}^5  \nonumber
\end{align}
Note that upon expressing a product as a linear combination of algebraic integers, which algebraic
degree over $\mathbf{Q}$ is $2^b$, then the integer coefficients of such linear combination sum to $2^{N-1-b}$. \\

Above we have shown the products with $q$-index $0$ (here omitted). The other $a$-indices for the products
$s_{0,2^b}^{b+1} \, s_{0,0}^{b+1}$ are obtained by cyclic shifts of the $a$-indices modulo $2^{b+1}$ on the lhs  and modulo $2^b$ on the rhs. The other $q$-indices in the products are also obtained by shifts, but not cyclic shifts. So on the right side, $q$-indices $\ge Q$ appear. Because the $q$-indices are defined in \ref{sqab_definition} as $0 \le q \le Q - 1$, these indices have to be reduced modulo $Q$. Reducing a $q$-index by $Q$ needs to increase the $a$-index by $+1$ to get the same sum. See
$\dots c({\color{red}q} + ({\color{red}a} + k \, 2^b) \, {\color{red}Q}) \dots$ in \ref{sqab_definition}. \\

Here an example. With a $q$-shift of $+4$, the blue indices have to be reduced. The $+1$-shifted $a$-values are shown in red.
\begin{align}
   & s_{4,4}^3 s_{4,0}^3 &    = \enspace  & (s_{4,2}^2) + (s_{\qsh[5],0}^2 + 2 s_{\qsh[5],2}^2 + s_{\qsh[5],3}^3) + (3 s_{\qsh[6],0}^3) \nonumber \\
   & s_{4,4}^3 s_{4,0}^3 &    = \enspace  & (s_{4,2}^2) + (s_{\qsh[0],0 \ash}^2 + 2 s_{\qsh[0],2 \ash}^2 + s_{\qsh[0],3 \ash}^3) + (3 s_{\qsh[1],0 \ash}^3) \nonumber
\end{align} \\

The $q$-shifts for the first product $s_1^1 s^1$ have a very simple form. Because the $a$-index on the rhs has only the value $0$, additional $a$-shifts with $+1$ yield no changes modulo $2^b = 1$. So these products are obtained by cyclic $q$-shifts. \\

{ \large \bf Explicit successive calculations of roots of quadratic polynomials} \\ 

The signs used to determine the signs of the square roots in the recursion now depend not only on the upper $b$-index but also on the lower $q$-index too: 
\begin{equation} \begin{array}{lllll}
   \tau^0 = \texttt{+}           & \tau_{1}^0 = \texttt{+}           & \tau_{2}^0 = \texttt{+}           & \tau_{3}^0 = \texttt{+}           & \tau_{4}^0 = \texttt{+} \\
   \tau^1 = \texttt{--}          & \tau_{1}^1 = \texttt{+-}          & \tau_{2}^1 = \texttt{--}          & \tau_{3}^1 = \texttt{++}          & \tau_{4}^1 = \texttt{-+} \\
   \tau^2 = \texttt{-++-}        & \tau_{1}^2 = \texttt{-+-+}        & \tau_{2}^2 = \texttt{----}        & \tau_{3}^2 = \texttt{---+}        & \tau_{4}^2 = \texttt{+--+} \\
   \tau^3 = \texttt{----\,-++-}  & \tau_{1}^3 = \texttt{-+--\,-++-}  & \tau_{2}^3 = \texttt{+--+\,++--}  & \tau_{3}^3 = \texttt{++--\,---+}  & \tau_{4}^3 = \texttt{-++-\,+-+-}
\end{array}  \nonumber \end{equation}
\begin{equation} \begin{array}{lll}
   \tau^4 = \texttt{+--+\,-++-\,----\,-+++}  & \tau_{1}^4 = \texttt{++-+\,---+\,+---\,-++-}  & \tau_{2}^4 = \texttt{+--+\,++++\,-+-+\,-+++} \\
	\tau_{3}^4 = \texttt{+-++\,+-++\,-++-\,---+}  & \tau_{4}^4 = \texttt{---+\,-++-\,-+-+\,-+-+}
\end{array}  \nonumber \end{equation}
\begin{equation} \begin{array}{ll}
   \tau^5 = \texttt{++++\,+---\,++-+\,+--+\,--+-\,--++\,--++\,-++-}     & \tau_{1}^5 = \texttt{+++-\,+++-\,---+\,++-+\,----\,++-+\,--++\,-+++}  \\
	\tau_{2}^5 = \texttt{+++-\,-++-\,+-+-\,++--\,-+++\,-+-+\,-+++\,-++-} & \tau_{3}^5 = \texttt{++++\,-++-\,+---\,----\,-+-+\,+-+-\,+--+\,+---}  \\
	\tau_{4}^5 = \texttt{+-+-\,++-+\,+++-\,++++\,+---\,--+-\,--+-\,+---}
\end{array}  \nonumber \end{equation}

\newpage
\noindent \textbf{\large Appendices}

\appendix

\section{Example: $P = 13 = 2^2 \cdot 3 + 1$, one trisection \label{examples_begin}}

$s_{0,0}^0,s_{1,0}^0, s_{2,0}^0$ are defined as roots of the cubic $s^3 + s^2 - 4 s + 1$.
Here, $2$ generates the multiplicative group. With trisection we get for these sums:
\begin{align}
	& C = \sqrt[3]{\frac{-5 + 3 \sqrt{3} \, i}{2 \, \sqrt{13}}}  \quad \vert C \vert = 1 \nonumber \\
	& s_{q,0}^0 = (- 1 + \sqrt{13} \, (\epsilon^{2+q} + \cc)) \, / \, 3 \quad \epsilon = e^{2 \pi i / 3}, \quad q = 0, 1, 2  \nonumber
\end{align}

\begin{align}
   & s_1^1 s^1 &    = \enspace  & s_{2,0}^0 \nonumber
\end{align}

Above, there is only one product, so cyclic $q$-shifts apply. \\

The signs for the single level (level $1$) are:
\begin{equation} \begin{array}{lll}
   \tau^0 = \texttt{+}           & \tau_{1}^0 = \texttt{+}           & \tau_{2}^0 = \texttt{+}
\end{array} \nonumber \end{equation}

\section{Example: $P = 29 = 2^2 \cdot 7 + 1$, three $7$-sections}

$s_{0,0}^0, \dots, s_{6,0}^0$ are defined as roots of the septic $s^7 + s^6 - 12 s^5 - 7 s^4 + 28 s^3 + 14 s^2 - 9 s + 1$.
Here, $2$ generates the multiplicative group. With 7-section, we get for these sums:
\begin{align}
	& \eta_n = 2 \, \cos (2 \pi \, 2^n / 7) \quad \sigma_n = \text{sign} (\sin (2 \pi \, 2^n / 7))\quad \sigma = [ +1, +1, -1 ] \quad n = 0, 1, 2 \nonumber \\
	& C_n = \sqrt[7]{\frac{5\,199 - 2\,597 \eta_n - 5\,831 \eta_n^2 - \sigma_n \, (4\,319 + 700 \eta_n - 1\,029 \eta_n^2) \, \sqrt{4 - \eta_n^2} \, i}{2 \cdot 29^{5/2}}} \quad \vert C_n \vert = 1 \nonumber \\
	& s_{q,0}^0 = (- 1 + \sqrt{29} \, (\epsilon^{1-q} \, C_0 + \epsilon^{3-2q} \, C_1 + \epsilon^{6-4q} \, C_2 + \cc)) \, /  \, 7 \quad \epsilon = e^{2 \pi i / 7}, \quad q = 0, \dots,  6  \nonumber
\end{align}
For the 7-section above, see also \cite{Adlaj_29}.

\begin{align}
   & s_1^1 s^1 &    = \enspace  & s_{4,0}^0 \nonumber
\end{align}

Above, there is only one product, so cyclic $q$-shifts apply. \\

The signs for the single level are:
\begin{equation} \begin{array}{lllllll}
   \tau^0 = \texttt{+}       & \tau_{1}^0 = \texttt{+}       & \tau_{2}^0 = \texttt{+}  & \tau_{3}^0 = \texttt{+}      & \tau_{4}^0 = \texttt{-}    & \tau_{5}^0 = \texttt{+}   & \tau_{6}^0 = \texttt{+} \\
\end{array} \nonumber \end{equation}

\section{Example: $P = 41 = 2^2 \cdot 5 + 1$, two quintisections}

$s_{0,0}^0, \dots, s_{4,0}^0$ are defined as roots of the quintic $s^5 + s^4 - 16 s^3 + 5 s^2 + 21 s - 9$.
Here, $6$ generates the multiplicative group. With quintisection, we get for these sums:
\begin{align}
	& C_{\pm} = \sqrt[5]{\frac{1\,962 \pm 50 \, \sqrt{5} + (255 \mp 75 \, \sqrt{5})) \, \sqrt{10 \pm 2 \, \sqrt{5}} \, i} {8 \cdot 41 ^ {3/2}}}  \quad \vert C_{\pm} \vert = 1 \nonumber \\
	& s_{q,0}^0 = (- 1 + \sqrt{41} \, ( \, \epsilon^{1-q} \, C_+
	                                    + \epsilon^{-2q} \, C_- + \cc \, )) \, / \, 5 \quad \epsilon = e^{2 \pi i / 5}, \quad q = 0, \dots,  4  \nonumber
\end{align}
For the quintisection above, see also \cite{Adlaj_41}.

\begin{align}
   & s_1^1 s^1 &    = \enspace  & s_{1,0}^0 + s_{2,0}^0 \nonumber \\
   & s_2^2 s^2 &    = \enspace  & s_{3,1}^1 \nonumber
\end{align}

$q$-shifts with additional $a$-shifts by $+1$ in red:
\begin{align}
	& q\text{-shift} \quad q \ge 2 & s_{q,2}^2 s_{q,0}^2 = s_{3+q\qsh[\, - \, Q],1 \ash}^1 \nonumber
\end{align}

The signs for the two levels are:
\begin{equation} \begin{array}{lllll}
   \tau^0 = \texttt{+}           & \tau_{1}^0 = \texttt{+}           & \tau_{2}^0 = \texttt{+}           & \tau_{3}^0 = \texttt{-}           & \tau_{4}^0 = \texttt{-} \\
   \tau^1 = \texttt{+-}          & \tau_{1}^1 = \texttt{++}          & \tau_{2}^1 = \texttt{-+}          & \tau_{3}^1 = \texttt{+-}          & \tau_{4}^1 = \texttt{+-}
\end{array}  \nonumber \end{equation}

\section{Example: $P = 97 = 2^5 \cdot 3 + 1$, one trisection \label{examples_end}}

$s_{0,0}^0,s_{1,0}^0, s_{2,0}^0$ are defined as roots of the cubic $s^3 + s^2 - 32s - 79$.
Here, $5$ generates the multiplicative group.  With trisection we get for these sums:
\begin{align}
	& C = \sqrt[3]{\frac{19 + 3 \sqrt{3} \, i}{2 \, \sqrt{97}}}  \quad \vert C \vert = 1 \nonumber \\
	& s_{q,0}^0 = (- 1 + \sqrt{97} \, (\epsilon^{1+q} \, C + \cc)) \, / \, 3 \quad \epsilon = e^{2 \pi i / 3}, \quad q = 0, 1, 2  \nonumber
\end{align}

\begin{align}
   & s_1^1 s^1 &    = \enspace  & 4 s^0 + 3 s_{1,0}^0 + s_{2,0}^0 \nonumber \\
   & s_2^2 s^2 &    = \enspace  & (s_{1}^1) + (s_{1,0}^1) + (s_{2,0}^1 + s_{2,1}^1) \nonumber \\
   & s_4^3 s^3 &    = \enspace  & (s_{1}^2) + (s_{2,0}^2) \nonumber \\
   & s_8^4 s^4 &    = \enspace  & s_{2,1}^3 \nonumber
\end{align}

$q$-shifts with additional $a$-shifts by $+1$ in red:
\begin{align}
	& \dots \nonumber \\
	& q\text{-shift} \quad q \ge 1  & s_{q,8}^4 s_{q,0}^4 = s_{2+q\qsh[\, - \, Q],1 \ash}^3 \nonumber
\end{align}

The signs for the four levels are:
\begin{equation} \begin{array}{lll}
   \tau^0 = \texttt{+}           & \tau_{1}^0 = \texttt{-}           & \tau_{2}^0 = \texttt{-} \\
   \tau^1 = \texttt{--}          & \tau_{1}^1 = \texttt{+-}          & \tau_{2}^1 = \texttt{-+} \\
   \tau^2 = \texttt{++++}        & \tau_{1}^2 = \texttt{+-+-}        & \tau_{2}^2 = \texttt{---+} \\
   \tau^3 = \texttt{+++-\,+--+}  & \tau_{1}^3 = \texttt{+---\,+--+}  & \tau_{2}^3 = \texttt{+++-\,-+++}
\end{array}  \nonumber \end{equation}

\section{Examples: $P = F_2, F_3,F_4$, the Fermat primes \label{special_fermat}}

\subsection{Example: $P = F_2 = 17 = 2^4 \cdot {\color{red} 1} + 1$ \label{special_17_Q1}}

This is a special case with $Q = {\color{red} 1}$ that allows us to construct the Fermat prime cases.
$s^0 = - 1$ is defined as a root of the degree $Q = 1$ polynomial $s + 1$.
Here $3$ generates the multiplicative group.

\begin{align}
   & s_1^1 s^1 &    = \enspace  & 4 s^0 \nonumber \\
   & s_2^2 s^2 &    = \enspace  & s^1 + s_1^1 \nonumber \\
   & s_4^3 s^3 &    = \enspace  & s_1^2 \nonumber
\end{align}

The signs for the three levels are:
\begin{equation} \begin{array}{lll}
   \tau^0 = \texttt{+}   & \tau^1 = \texttt{++}     & \tau^2 = \texttt{++--}
\end{array} \nonumber \end{equation}

\subsection{Example: $P = F_2 = 17 = 2^3 \cdot {\color{red} 2} + 1$, one trivial bisection \label{special_17_Q2}}

This is a special case with $Q = {\color{red} 2}$ that allows us to construct the Fermat prime cases. The bisection is trivial: a zero angle is to be bisected.
$s_{0,0}^0,s_{1,0}^0$ are defined as roots of the quadratic $s^2 + s - 4$.
Here, $3$ generates the multiplicative group. With trivial bisection we get for these sums:
\begin{align}
	& C = \sqrt[2]{1}  \quad \vert C \vert = 1 \nonumber \\
	& s_{q,0}^0 = (- 1 + \sqrt{17} \, \epsilon^{q} \, C) \, / \, 2 \quad \epsilon = e^{2 \pi i / 2} = -1, \quad q = 0, 1  \nonumber
\end{align}

\begin{align}
   & s_1^1 s^1 &    = \enspace  & s^0 + s_{1,0}^0 \nonumber \\
   & s_2^2 s^2 &    = \enspace  & s_{1,0}^1 \nonumber
\end{align}

$q$-shifts with an additional $a$-shift by $+1$ in red:
\begin{align}
	& q\text{-shift} \quad q = 1  & s_{1,1}^3 s_{1,0}^3 = s_{\qsh[0],0 \ash}^2 \nonumber
\end{align}

The signs for the two levels are:
\begin{equation} \begin{array}{llll}
   \tau^0 = \texttt{+}   & \tau_{1}^0 = \texttt{+}     & \tau^1 = \texttt{+-}   & \tau_{1}^1 = \texttt{+-}
\end{array} \nonumber \end{equation}

\subsection{Example: $P = F_2 = 17 = 2^2 \cdot {\color{red} 4} + 1$, one  bisection \label{special_17_Q4}}

This is a special case with $Q = {\color{red} 4}$ that allows to construct the Fermat prime cases.
$s_{0,0}^0, \dots,  s_{3,0}^0$ are defined as roots of the quartic $s^4 + s^3 - 6 s^2 - s +1$.
Here, $3$ generates the multiplicative group.
\begin{align}
	& C = \sqrt[4]{\frac{(1 - 4 \, i)^2}{17}}  \quad \vert C \vert = 1 \quad \text{this is an angle bisection}\nonumber \\
	& s_{q,0}^0 = (- 1 + \sqrt{17} \, ( (\epsilon^{-q} \, C + \cc) + \epsilon^{-2q} )) \, / \, 4 \quad \epsilon = e^{2 \pi i / 4} = +i, \quad q = 0, \dots, 3  \nonumber
\end{align}

\begin{align}
   & s_1^1 s^1 &    = \enspace  & s_{1,0}^0 \nonumber
\end{align}

Above only one product, so cyclic $q$-shifts. \\

The signs for the single level are:
\begin{equation} \begin{array}{llll}
   \tau^0 = \texttt{+}   & \tau_{1}^0 = \texttt{+}     & \tau_{2}^0 = \texttt{-}   & \tau_{3}^0 = \texttt{-}
\end{array} \nonumber \end{equation}

\subsection{Example: $P = F_2 = 17 = 2^1 \cdot {\color{red} 8} + 1$, two quartisections \label{special_17_Q8}}

This is a special case with $Q = {\color{red} 8}$ that allows to construct the Fermat prime cases.
$s_{0,0}^0, \dots,  s_{7,0}^0$ are defined as roots of the octic $s^8 + s^7 - 7 s^6 - 6 s^5 + 15 s^4 + 10 s^3 - 10 s^2 - 4 s + 1$.
Here, $3$ generates the multiplicative group.
\begin{align}
	& C_{\pm} = \sqrt[8]{\frac{(- 1 \mp 48 \, \sqrt{2} + ( \pm 4 - 12 \, \sqrt{2}) \, i)^2}{17^3}}  \quad C_2 = \sqrt[\color{red} \bf 4]{\frac{(1 - 2^2 \, i)^2}{17}}  \quad \vert C_{\pm} \vert, \vert C_2 \vert  = 1 \nonumber \\
	& s_{q,0}^0 = (- 1 + \sqrt{17} \, ( (\epsilon^{7-q} \, C_+ + \epsilon^{2-2q} \, C_2 + \epsilon^{2-3q} \, C_- +  \cc) + \epsilon^{-4q} )) \, / \, 8 \nonumber \\
	&  \qquad \quad \epsilon = e^{2 \pi i / 8} = (1 + i) / \sqrt{2}, \quad q = 0, \dots, 7  \nonumber
\end{align}

There are no more quadratic extensions.

\subsection{Example: $P = F_3 = 257 = 2^5 \cdot {\color{red} 8} + 1$, two quartisections \label{special_257_Q8}}

This is a special case with $Q = {\color{red} 8}$, that allows to construct the Fermat prime cases.
$s_{0,0}^0, \dots,  s_{7,0}^0$ are defined as roots of the octic $s^8 + s^7 - 112 s^6 - 542 s^5 + 349 s^4 + 4481 s^3 + 778 s^2 - 9316 s + 2312$.
Here $3$ generates the multiplicative group.
\begin{align}
	& C_{\pm} = \sqrt[8]{\frac{(193 \mp  1\,920 \, \sqrt{2} - ( \pm 3\,088 + 120 \, \sqrt{2}) \, i)^2}{P^3}}  \quad C_2 = \sqrt[\color{red} \bf 4]{\frac{(1 - 2^4 \, i)^2}{P}}  \quad \vert C_{\pm} \vert, \vert C_2 \vert  = 1 \nonumber \\
	& s_{q,0}^0 = (- 1 + \sqrt{P} \, ( (\epsilon^{-q} \, C_+ + \epsilon^{2-2q} \, C_2 + \epsilon^{7-3q} \, C_- +  \cc) + \epsilon^{-4q} )) \, / \, 8 \nonumber \\
	&  \qquad \qquad \qquad \epsilon = e^{2 \pi i / 8} = (1 + i) / \sqrt{2}, \quad q = 0, \dots, 7  \nonumber
\end{align}

\begin{align}
   & s_1^1 s^1 &    = \enspace  & 2 s^0 + 2 s_{2,0}^0 + s_{4,0}^0 + 2 s_{5,0}^0 + s_{6,0}^0 \nonumber \\
   & s_2^2 s^2 &    = \enspace  & s^1 + s_{1,0}^1 + s_{2,0}^1 + s_{5,0}^1 \nonumber \\
   & s_4^3 s^3 &    = \enspace  & s_{1,0}^2 + s_{7,2}^2 \nonumber \\
   & s_8^4 s^4 &    = \enspace  & s_{0,7}^3 \nonumber
\end{align}

The signs for the four levels are: \dots

\subsection{Example: $P = F_4 = 65\,537 = 2^{13} \cdot {\color{red} 8} + 1$, two quartisections \label{special_257_Q8}}

This is a special case with $Q = {\color{red} 8}$ that allows to construct the Fermat prime cases.
$s_{0,0}^0, \dots,  s_{7,0}^0$ are defined as roots of the octic $s^8 + s^7 - 28672 s^6 - 2104352 s^5 + 23067904 s^4 + 4319150336 s^3 + 4832370688 s^2 - 2217310683136 s + 18142201905152$.
Here, $3$ generates the multiplicative group.
\begin{align}
	& C_{\pm} = \sqrt[8]{\frac{(64\,513 \mp  2\,088\,960 \, \sqrt{2} - ( \pm 16\,515\,328 + 8\,160 \, \sqrt{2}) \, i)^2}{P^3}}  \nonumber \\
	&  C_2 = \sqrt[\color{red} \bf 4]{\frac{(1 - 2^8 \, i)^2}{P}}  \quad \vert C_{\pm} \vert, \vert C_2 \vert  = 1 \nonumber \\
	& s_{q,0}^0 = (- 1 + \sqrt{P} \, ( (\epsilon^{2-q} \, C_+ + \epsilon^{6-2q} \, C_2 + \epsilon^{5-3q} \, C_- +  \cc) + \epsilon^{-4q} )) \, / \, 8 \nonumber \\
	&  \qquad \qquad \qquad \epsilon = e^{2 \pi i / 8} = (1 + i) / \sqrt{2}, \quad q = 0, \dots, 7  \nonumber
\end{align}

\begin{align}
   & s_1^1 s^1 &    = \enspace  & 284 s^0 + 237 s_{1,0}^0 + 272 s_{2,0}^0 + 237 s_{3,0}^0 + 256 s_{4,0}^0 + 269 s_{5,0}^0 \nonumber \\
	&           &                & + 256 s_{6,0}^0 + 237 s_{7,0}^0 \nonumber \\
   &           &    \dots \nonumber \\
   & s_{1024}^{11} s^{11} &    = \enspace  & s_{4,395}^{10} + s_{4,1012}^{10} \nonumber \\
   & s_{2048}^{12} s^{12} &    = \enspace  & s_{0,384}^{11} \nonumber
\end{align}

The signs for the twelve levels are: \dots

\section{Examples: $P = 13,97$ with $Q = 6, 12$.}

\subsection{Example: $P = 13 = 2^1 \cdot {\color{red} 6} + 1$, two trisections \label{special_13_Q6}}

This is a special case with $Q = {\color{red} 6}$.
$s_{0,0}^0, \dots,  s_{5,0}^0$ are defined as roots of the sextic $s^6 + s^5 - 5 s^4 - 4 s^3 + 6 s^2 + 3 s - 1$.
Here, $2$ generates the multiplicative group.
\begin{align}
	& C_1 = \sqrt[\color{red} \bf 3]{\frac{1 - 15 \, \sqrt{3} \, i}{2 \cdot 13}}  \quad C_2 = \sqrt[6]{\frac{(5 + 3 \, \sqrt{3} \, i)^2}{4 \cdot 13}}  \quad \vert C_1 \vert, \vert C_2 \vert  = 1 \nonumber \\
	& s_{q,0}^0 = (- 1 + \sqrt{13} \, ( (\epsilon^{1-q} \, C_1 + \epsilon^{1-2q} \, C_2 +  \cc) + \epsilon^{-3q} )) \, / \, 6 \nonumber \\
	&  \qquad \quad \epsilon = e^{2 \pi i / 6} = (1 + \sqrt{3} \, i) / 2, \quad q = 0, \dots, 5  \nonumber
\end{align}

There are no more quadratic extensions.

\subsection{Example: $P = 97 = 2^4 \cdot {\color{red} 6} + 1$, two trisections \label{special_97_Q6}}

This is a special case with $Q = {\color{red} 6}$.
$s_{0,0}^0, \dots,  s_{5,0}^0$ are defined as roots of the sextic $s^6 + s^5 - 40 s^4 + 45 s^3 + 236 s^2 - 230 s - 389$.
Here $5$ generates the multiplicative group.
\begin{align}
	& C_1 = \sqrt[\color{red} \bf 3]{\frac{167 - 57 \, \sqrt{3} \, i}{2 \cdot 97}}  \quad C_2 = \sqrt[6]{\frac{(19 - 3 \, \sqrt{3} \, i)^2}{4 \cdot 97}}  \quad \vert C_1 \vert, \vert C_2 \vert  = 1 \nonumber \\
	& s_{q,0}^0 = (- 1 + \sqrt{97} \, ( (\epsilon^{5-q} \, C_1 + \epsilon^{3-2q} \, C_2 +  \cc) + \epsilon^{-3q} )) \, / \, 6 \nonumber \\
	&  \qquad \quad \epsilon = e^{2 \pi i / 6} = (1 + \sqrt{3} \, i) / 2, \quad q = 0, \dots, 5  \nonumber
\end{align}

\begin{align}
   & s_1^1 s^1 &    = \enspace  & s_{1,0}^0 + s_{2,0}^0 + s_{3,0}^0 + s_{5,0}^0 \nonumber \\
   & s_2^2 s^2 &    = \enspace  & s_{2,0}^1 + s_{3,0}^1 \nonumber \\
   & s_4^3 s^3 &    = \enspace  & s_{5,0}^2 \nonumber
\end{align}

The signs for the $3$ levels are:
\begin{equation} \begin{array}{lll}
   \tau^0     = \texttt{-}       & \tau_{1}^0 = \texttt{+}           & \tau_{2}^0 = \texttt{-} \\
   \tau_{3}^0 = \texttt{-}       & \tau_{4}^0 = \texttt{-}           & \tau_{5}^0 = \texttt{+} \\
   \tau^1     = \texttt{++}      & \tau_{1}^1 = \texttt{++}          & \tau_{2}^1 = \texttt{--} \\
   \tau_{3}^1 = \texttt{++}      & \tau_{4}^1 = \texttt{--}          & \tau_{5}^1 = \texttt{-+} \\
   \tau^2     = \texttt{+++-}    & \tau_{1}^2 = \texttt{+-+-}        & \tau_{2}^2 = \texttt{++-+} \\
   \tau_{3}^2 = \texttt{+--+}    & \tau_{4}^2 = \texttt{---+}        & \tau_{5}^2 = \texttt{+-++}
\end{array}  \nonumber \end{equation}

\subsection{Example: $P = 97 = 2^3 \cdot {\color{red} 12} + 1$, three sextisections \label{special_97_Q12}}

This is a special case with $Q = {\color{red} 12}$.
$s_{0,0}^0, \dots,  s_{11,0}^0$ are defined as roots of the degree $12$ polynomial $s^{12} + s^{11} - 44 s^{10} - 23 s^9 + 608 s^8 + 288 s^7 - 3367 s^6 - 1647 s^5 + 7459 s^4 + 2633 s^3 - 7037 s^2 - 1034 s + 2209$.
Here, $5$ generates the multiplicative group.
\begin{align}
	& C_{\pm} = \sqrt[12]{\frac{(49\,599 \pm 51\,756 \, \sqrt{3} + (151\,636 \mp 16\,929  \, \sqrt{3}) \, i)^2}{4 \cdot 97^5}} \nonumber \\
	& C_4 = \sqrt[\color{red} \bf 6]{\frac{167 - 57 \, \sqrt{3} \, i}{2 \cdot 97}} \nonumber \quad C_3 = \sqrt[\color{red} \bf 4]{\frac{(9 + 4 \, i)^2}{97}} \quad \vert C_{\pm} \vert, \vert C_4 \vert, \vert C_3 \vert  = 1 \nonumber \\
	& s_{q,0}^0 = (- 1 + \sqrt{97} \, ( (\epsilon^{8-q} \, C_+ + \epsilon^{-2q} \, C_4^2 + \epsilon^{3-3q} \, C_3 + \epsilon^{8-4q} \, C_4  \nonumber \\
	& \quad \qquad  + \epsilon^{5-5q} \, C_- +  \cc) + \epsilon^{-6q} )) \, / \, 12 \nonumber \\
	&  \qquad \epsilon = e^{2 \pi i / 12} = (\sqrt{3} + i) / 2, \quad q = 0, \dots, 11  \nonumber
\end{align}

\begin{align}
   & s_1^1 s^1 &    = \enspace  & s_{2,0}^0 + s_{3,0}^0 \nonumber \\
   & s_2^2 s^2 &    = \enspace  & s_{5,0}^1 \nonumber
\end{align}

The signs for the $2$ levels are:
\begin{equation} \begin{array}{llllll}
     \tau^0     = \texttt{+}       & \tau_{1}^0 = \texttt{+}           & \tau_{2}^0 = \texttt{-}
   & \tau_{3}^0 = \texttt{+}       & \tau_{4}^0 = \texttt{-}           & \tau_{5}^0 = \texttt{-} \\
     \tau_{6}^0 = \texttt{+}       & \tau_{7}^0 = \texttt{+}           & \tau_{8}^0 = \texttt{-}
   & \tau_{9}^0 = \texttt{+}       & \tau_{10}^0 = \texttt{-}          & \tau_{11}^0 = \texttt{+} \\
     \tau^1     = \texttt{++}      & \tau_{1}^1 = \texttt{++}          & \tau_{2}^1 = \texttt{--}
   & \tau_{3}^1 = \texttt{+-}      & \tau_{4}^1 = \texttt{--}          & \tau_{5}^1 = \texttt{++} \\
     \tau_{6}^1 = \texttt{+-}      & \tau_{7}^1 = \texttt{--}          & \tau_{8}^1 = \texttt{++}
   & \tau_{9}^1 = \texttt{-+}      & \tau_{10}^1 = \texttt{-+}         & \tau_{11}^1 = \texttt{-+}
\end{array}  \nonumber \end{equation}

\vspace{1.0cm}

\bibliographystyle{amsplain}

\end{document}